\documentclass[a4paper,10pt]{article}
\usepackage{subfigure}
\usepackage{graphicx}
\usepackage{epstopdf}
\usepackage{multirow}
\usepackage{amsmath}
\usepackage{amsfonts}
\usepackage{amssymb}
\usepackage{amsthm}
\usepackage{marvosym}
\usepackage{tikz}
\usetikzlibrary{arrows,topaths}
\usepackage{tkz-berge}
\usepackage{tkz-graph}
\usepackage{enumerate}
\usepackage{indentfirst}
\usepackage{textcomp}
\usepackage{latexsym}
\usepackage{amssymb}
\usepackage{amsmath}
\usepackage{algorithm}
\usepackage{algorithmic}

\newtheorem{defn}{Definition}
\newtheorem{thm}[defn]{Theorem}

\begin{document}
\date{}
\title{\bf\large{A note on ¡°Graphical Notation Reveals Topological Stability Criteria for Collective Dynamics in Complex Network¡±}
\footnote{Partially supported by a grant from China Scholarship Council and National Natural Science Foundation of China (11101071, 11271001) and Fundamental Research Funds for the Central Universities (ZYGX2016J131, ZYGX2016J138).}}
\author{Peng-hui He$^1$, Hou-biao Li$^1$\footnote{Corresponding Author, Email: lihoubiao0189@163.com}, Hong Li$^1$, Nan Jiang$^2$
 \\
{\small $^1$School of Mathematical Sciences, University of Electronic Science}\\
{\small and Technology of China, Chengdu, 611731, P.R. China.}\\
{\small $^2$School of Computer Science and Engineering, University of Electronic Science}\\
{\small and Technology of China, Chengdu, P.R. China}}

\maketitle

\begin{abstract}
This paper clarifies the main research methods and ideas of the thesis \cite{1,2,3}. The special calculation process is also realized by corresponding computer algorithm. Finally, we introduce zero rows sum case and give the corresponding algorithm, which greatly simplifies the relevant calculation of paper \cite{1,3}.
\end{abstract}

\section{Introduction}
To reveal topological stability criteria for collective dynamics in complex network, the paper \cite{1} proposed a Coates graphical notation by which certain spectral properties of complex systems can be
rewritten concisely and interpreted topologically. For deriving topological stability criteria, the Jacobi¡¯s signature criterion (JSC), also called Sylvester criterion is applied, which states that a Hermitian matrix $J$ with rank $r$ is negative definite if and only if all principal minors of order $q \le r$  have the sign of ${(-1)^q}$, i.e.,
\begin{equation}\label{eq:1.1}
{\mathop{\rm sgn}} ({D_{\left| S \right|,S}}) = {( - 1)^{\left|S \right|}},{\rm{  }}\forall S\subset \{ 1,2, \cdots,n\} ,\left| S \right| \le r,
\end{equation}
where ${D_{\left| S \right|,S}}:=det(J_{ik}), i,k\in S$.
Therefore, the evaluation of determinants of Jacobian matrices $J\in \mathbb{R}^{n\times n}(J_{ik}=\partial \dot{x}_i/\partial \dot{x}_k)$ is important and essential.

For systems with zero row sum (ZRS), the
minors (\ref{eq:1.1}) can be expresses as
\begin{equation}\label{eq:1.2}
{D_{\left| S \right|,S}}=(-1)^{\left|S \right|}\Phi_S,
\end{equation}
where $\Phi_S$ is defined as the sum over all
products corresponding to all acyclic graphs that can be drawn by placing $\left|S\right|$ (undirected) links such that each link starts at a different node in $S$ in Ref. \cite{1}.

In this note, we will clarify the main idea of the $\Phi_S$ and propose its corresponding algorithm and perfect the corresponding theory.

The remainder of the paper is organized as follows. In Section 2, we mainly exhibit the precise definition of the $\Phi_S$ by a specific example and give its corresponding algorithm to complement the corresponding results. In addition, a more simpler result for general matrix is introduced and the related topological algorithm is given in Section 3. Finally, some concluding remarks on further research are given in Section 4.



\section{Definition of the $\Phi_S$ in symmetric case}

First, let us recall the corresponding results of paper \cite{1,2}.
Matrix $A$ is a symmetric with zero row sum $n \times n$ matrix.
$G= {\cal G}(A)$ is the Coates graph of $A$. In paper \cite{2}, the $\Phi_S$ is defined via the three criteria (I), (II) and (III):
\begin{description}
  \item[(I)] the sub-graph contains exactly $\left|S\right|$ links;
  \item[(II)] there is at least one link connecting to every node in $S$;
  \item[(III)] every link connects to at least one node in $S$.
\end{description}

In paper \cite{1}, the above three criteria is substituted by the following definition:
''define $\Phi_S$ as the sum over all
products corresponding to all acyclic graphs that can be drawn by placing $\left|S\right|$ (undirected) links such that each link starts at a different node in $S$.''
Note that the condition is very important
\begin{description}
  \item[(IV)] each link starts at a different node in $S$, which excludes the extra subgraphs (see Figure 15).
\end{description}

Recently, another similar definition is stated in Ref. \cite{22}: For $S \subseteq \{ 1,2,3, ... ,n\}$, $S \ne \emptyset$,
$$
{{\cal F}_S}: = \left\{ {K \subseteq E(G)\left| {|K| = |S|} \right.} \right\}$$
where each connected component in $K$ contains at least one vertex not in $S$. In fact, what is $\Phi_S$ exactly definition? and how is it calculated for large matrices? Next, Let us explore this problem and give a feasible calculation plan.

Obviously, in practical computation, the above some conditions are not easy to be judged in graphs. Sometimes it is easier to translate it into a "matrix form", or "two-dimensional array" (see Algorithm 1). For example, the condition (IV) is actually equivalent to removing subgraphs whose row indices are in $S$ and the number of non-zero elements whose column indices are outside $S$ is greater than or equal to 2. For Example 1, Figure 15 corresponds to the following matrix.
\[\left[ {\begin{array}{*{20}{c}}
1&2&0& \vdots &0&0\\
2&3&0& \vdots &0&0\\
0&0&4& \vdots &5&{ - 8}\\
 \cdots & \cdots & \cdots & \cdots & \cdots & \cdots \\
0&0&5& \vdots &0&0\\
0&0&{ - 8}& \vdots &0&0
\end{array}} \right].\]

\begin{algorithm}[htb]
\caption{PRL symmetric algorithm}
\label{alg:1}
\begin{algorithmic}[1]
\REQUIRE
$A$, $S=\{ {n_1},{n_2},{n_3},...,{n_q}\} $.
\ENSURE
\STATE Create a undirected Coates graph of matrix $A$, Initialization ${[A]_{S,S}} = 0$;
\STATE Preprocessing-removes the edges not adjacent to any vertex in the $S$ collection;
\STATE Map the edge information ($e[i][j]$) of the preprocessed graph to a one-dimensional array $b[k]$;
\STATE With recursive implementation, from $B$ arbitrary selection of $\left|S\right|$ elements, and stored in the $a[]$;
\STATE Restore the $a[]$ to the two-dimensional array $e[i][j]$, and use it as the side information of the subgraph to construct the subgraph SubG;
\STATE Determine whether there is a loop in subgraph SubG with the depth-first search algorithm. If there is a ring, then exit; otherwise, judge the subgraph SubG whether all vertices of $S$ are connected with an edge and each link starts at a different node in $S$; If yes, then calculate ${[A]_{S,S}}={[A]_{S,S}}+\prod\limits_{e \in SubG}{\omega (e)} $, exit.
\RETURN ${[A]_{S,S}}$.
\end{algorithmic}
\end{algorithm}

Next, let us look at an example.

{\bfseries Example 1}. Let $A$ be a symmetric square matrix
\[A = \left[ {\begin{array}{*{20}{c}}
1&2&0&0&{ - 3}\\
2&3&{ - 1}&0&{ - 4}\\
0&{ - 1}&4&5&{ - 8}\\
0&0&5&6&{ - 11}\\
{ - 3}&{ - 4}&{ - 8}&{ - 11}&{26}
\end{array}} \right].\]
$S=\left\{{1,2,3}\right\}$. The Coates graph of $A$ is shown Figure 1.

\begin{figure}
\centering
\includegraphics[width=1.2in]{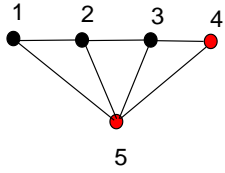}
\caption{The Coates graph of $A$}
\end{figure}

${\phi _{3,S}}$: According to Algorithm 1, all sub-graphs meet the conditions (I)-(III) have 14, which are shown as follows (see Figure 2-15). However, duo to the condition (IV), the Figure 15 is excluded. Therefore, according to (\ref{eq:1.2}), the determent ${\rm{D}}_{3,S}$ is
\[{{\rm{D}}_{3,S}}=(-1)^{3}(6+8+16+64-24-12+48-96-10-30-40+15+60)=-5.\]

\begin{figure}
\begin{minipage}[t]{0.5\linewidth}
\centering
\includegraphics[width=1.0in]{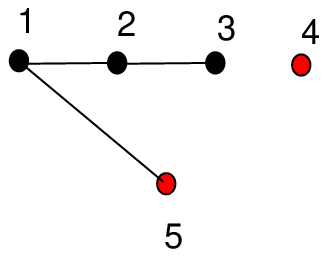}
\caption{$(-3)\times(2)\times (-1)=6 $}
\end{minipage}%
\begin{minipage}[t]{0.5\linewidth}
\centering
\includegraphics[width=1.0in]{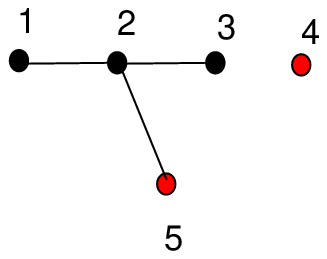}
\caption{$(2)\times(-1)\times (-4)=8$ }
\end{minipage}
\end{figure}

\begin{figure}
\begin{minipage}[t]{0.5\linewidth}
\centering
\includegraphics[width=1.0in]{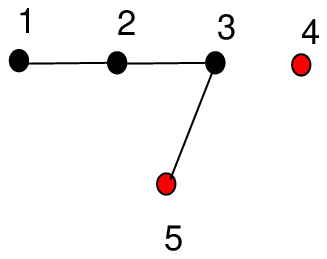}
\caption{$(2)\times(-1)\times (-8)=16 $}
\end{minipage}
\begin{minipage}[t]{0.5\linewidth}
\centering
\includegraphics[width=1.0in]{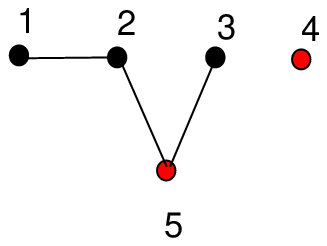}
\caption{$(2)\times(-4)\times (-8)=64$}
\end{minipage}%
\end{figure}

\begin{figure}
\begin{minipage}[t]{0.5\linewidth}
\centering
\includegraphics[width=1.0in]{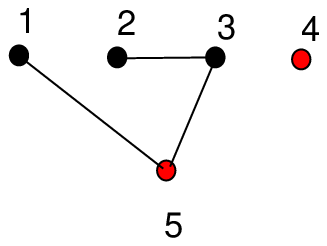}
\caption{$(-3)\times(-1)\times (-8)=-24$}
\end{minipage}
\begin{minipage}[t]{0.5\linewidth}
\centering
\includegraphics[width=1.0in]{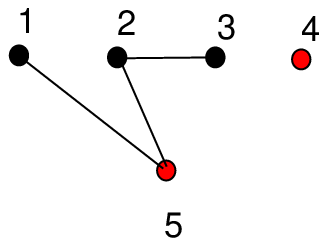}
\caption{$(-3)\times(-4)\times (-1)=-12$ }
\end{minipage}
\end{figure}

\begin{figure}
\begin{minipage}[t]{0.5\linewidth}
\centering
\includegraphics[width=1.0in]{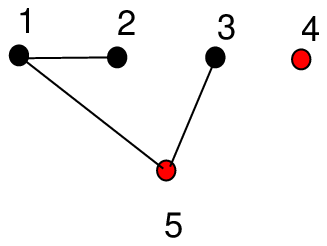}
\caption{$(-3)\times(2)\times (-8)=48$ }
\end{minipage}%
\begin{minipage}[t]{0.5\linewidth}
\centering
\includegraphics[width=1.0in]{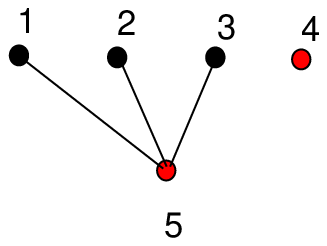}
\caption{$(-3)\times(-4)\times (-8)=-96 $}
\end{minipage}
\end{figure}

\begin{figure}
\begin{minipage}[t]{0.5\linewidth}
\centering
\includegraphics[width=1.0in]{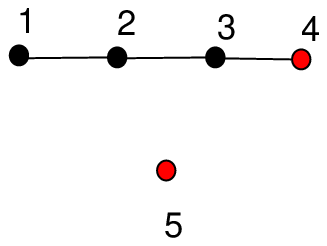}
\caption{$(5)\times(2)\times(-1)=-10$}
\end{minipage}
\begin{minipage}[t]{0.5\linewidth}
\centering
\includegraphics[width=1.0in]{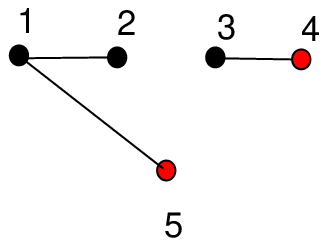}
\caption{$(-3)\times(2)\times(5)=-30$ }
\end{minipage}%
\end{figure}

\begin{figure}
\begin{minipage}[t]{0.5\linewidth}
\centering
\includegraphics[width=1.0in]{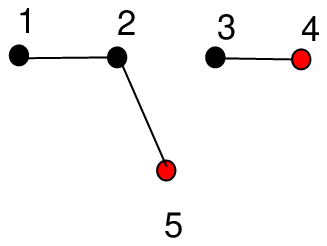}
\caption{$(-4)\times(2)\times(5)=40$}
\end{minipage}
\begin{minipage}[t]{0.5\linewidth}
\centering
\includegraphics[width=1.0in]{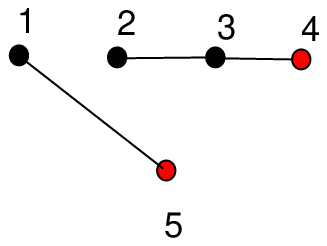}
\caption{$(-3)\times(5)\times(-1)=15$ }
\end{minipage}%
\end{figure}

\begin{figure}
\begin{minipage}[t]{0.5\linewidth}
\centering
\includegraphics[width=1.0in]{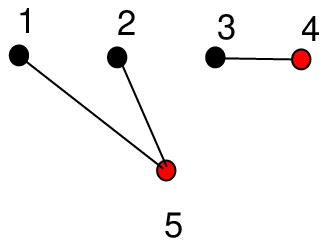}
\caption{$(-3)\times(-4)\times(5)=60$ }
\end{minipage}
\begin{minipage}[t]{0.5\linewidth}
\centering
\includegraphics[width=1.0in]{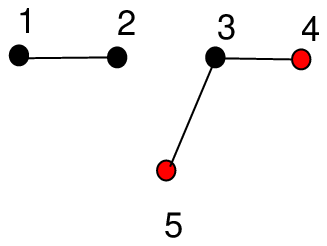}
\caption{$(2)\times(5)\times(-8)=-80$ }
\end{minipage}%
\end{figure}

\section{Definition of the $\Phi_S$ in asymmetric case}
Next, let us generalize the above result to general matrix. Main results come from the book \cite{3}. First, we introduce the concepts of directed tree with reference node $i$ and equicofactor matrix.

\begin{defn} [Directed tree \cite{3}] A subgraph, denoted by the symbol ${t_i}$, of a directed graph $G$ is called a directed tree of $G$ with reference node $i$ if and only if
\begin{description}
  \item[(1)] ${t_i}$ is a connected acyclic subgraph and contains all the nodes of $G$, i.e., it is a tree of $G$;
  \item[(2)] the outgoing degree of each node of ${t_i}$ is 1 except the node $i$ which has outgoing degree is 0.
\end{description}
\end{defn}

\begin{defn}
[Equicofactor matrix \cite{3}] A square matrix is said to be an equicofactor matrix if it has the properties that the sum of elements of every row and every column equals zero.
\end{defn}

If $A$ is an equicofactor matrix, then all the cofactors of the elements of $A$ are equal, and the following famous theorem holds.

\begin{thm}[Chen Wai-kai Theorem \cite{3}] Matrix $A$ is an equicofactor matrix, ${A_{i,j}}(i,j = 1,2, \cdots,n)$ is the first order cofactor of matrix $A$. then ${A_{i,j}}$ is equal to the sum of the weight values of all the directed trees with reference any node $i$, i.e.,
\[{A_{i,j}} = \sum\limits_{{t_k}} {f({t_k})} ,\]
where ${t_k}$ is $G$'s directed tree with reference node $k$.
\end{thm}

In the actual situation, how do we determine whether all directed trees with reference node $i$ have been fully searched? There exist some corresponding results in book \cite{3}. In addition, we can also obtain an accurate estimate according to Laplace matrix.

As it is well known, the Laplace matrix ${\cal L}\left( G \right)$ of graph $G$ is defined by ${\cal L}\left( G \right) = {\cal D}\left( G \right) - {\cal A}\left( G \right)$, where ${\cal D}\left( G \right)$ is the degree matrix and ${\cal A}\left( G \right)$ is the adjacency matrix of the graph $G$ without weight.

\begin{thm}[\cite{GodsilRoyle2001}] $G$ is the accompanying digraph of $A$, the number of all the directed trees of which taking any one node as its root is the same. The number is denoted by $N$:
\[N = \frac{1}{n}\prod\limits_i^q {{\lambda _i}} \]
where $n$ is the dimension of ${\cal L}(G)$, ${\lambda _i}(i = 1,2, \cdots,q)$ is the nonzero eigenvalues of ${\cal L}\left( G \right)$.
\end{thm}

Next, let us consider the general matrix. Set $A$ be an arbitrary square matrix.
\[A = \left[ {\begin{array}{*{20}{c}}
{{a_{11}}}&{{a_{12}}}& \cdots &{{a_{1n}}}\\
{{a_{21}}}&{{a_{22}}}& \cdots &{{a_{2n}}}\\
 \vdots & \vdots & \cdots & \vdots \\
{{a_{n1}}}&{{a_{n2}}}& \cdots &{{a_{nn}}}
\end{array}} \right]\]
In order to make use of the Chen's theorem, we extend $A$ to an $n+1$  order square matrix $\tilde{A}$.
\[\tilde A = \left[ {\begin{array}{*{20}{c}}
{{a_{11}}}&{{a_{12}}}& \cdots &{{a_{1n}}}&{ - \sum\limits_{j = 1}^n {{a_{1j}}} }\\
{{a_{21}}}&{{a_{22}}}& \cdots &{{a_{2n}}}&{ - \sum\limits_{j = 1}^n {{a_{2j}}} }\\
 \cdots & \cdots & \cdots & \cdots & \cdots \\
{{a_{n1}}}&{{a_{n2}}}& \vdots &{{a_{nn}}}&{ - \sum\limits_{j = 1}^n {{a_{nj}}} }\\
{ - \sum\limits_{j = 1}^n {{a_{j1}}} }&{ - \sum\limits_{j = 1}^n {{a_{j2}}} }& \cdots &{ - \sum\limits_{j = 1}^n {{a_{jn}}} }&{\sum\limits_{i = 1}^n {\sum\limits_{j = 1}^n {{a_{ij}}} } }
\end{array}} \right].\]
It is easy to get that matrix $\tilde{A}$ is an equicofactor matrix. Matrix $A$ is just the first order cofactor of matrix $\tilde{A}$. So, according to Chen's theorem,
$$\det(A)=\tilde{A}_{n+1,n+1}=\sum_{t_k} f(t_{k}) $$
Therefore, we transform the value problem of determinant of $A$ into looking for all the directed trees of which taking any one node as its root in the associated graph $G$ of $\tilde{A}$ .

So far, the best algorithm for finding all spanning trees in undirected and directed graphs is the one given by Gabow and Meers in 1978 (hereinafter referred to as G-M algorithm) \cite{4}, which is an $\mathcal{O}(NE+V+E)$ time complexity algorithm. Here $N$ refers to the number of spanning trees of a graph having $V$ vertices and $E$ edges. Recently, some improved versions have been presented in \cite{tree,44}, respectively. For example, the algorithm \cite{tree} takes $\mathcal{O}(NV+V^3)$ time and $O(V^2)$ space for enumerating all the spanning trees. But they is usually more complicated. For convenience, this paper uses G-M algorithm to find all spanning trees of digraph. In this way, we can get the topological algorithm of the asymmetric matrix, which is shown in Algorithm 2.

\begin{algorithm}[htb]
\caption{The Topological Algorithm of Determinant.}
\label{alg:2}
\begin{algorithmic}[1]
\REQUIRE
Input $n\times n$ matrix $A$.
\ENSURE
\STATE Create a matrix $\tilde{A}$ according to $A$ such that $det(A)=\tilde{A}_{n+1,n+1}$;
\STATE Create a directed Coates graph of matrix $\tilde{A}$;
\STATE According to the G-M algorithm, we get all the $G$ directed trees $T=\{T_{1},T_{2},...,T_{N}\}$ with reference certain node, where $T_{i}(i=1,2,...,N)$ is a directed tree, $N$ is the number of all the $G$ directed trees with reference this node;
\STATE Calculate $\det(A)=\sum_{T_{i}\in T}{\prod_{e\in{T_{i}}}{\omega (e)}}$;
\RETURN $det(A)$.
\end{algorithmic}
\end{algorithm}

Now, let us look back at Example 1. Obviously,
\[{D_{3,S}} = \left( {\begin{array}{*{20}{c}}
1&2&0\\
2&3&{ - 1}\\
0&{ - 1}&4
\end{array}} \right)\]
Equicofactor matrix $\tilde{D}_{3,S}$ is constructed by $D_{3,S}$ as follows.
\[\tilde{D}_{3,S}= \left[ {\begin{array}{*{20}{c}}
1&2&0&{ - 3}\\
2&3&{ - 1}&{ - 4}\\
0&{ - 1}&4&{ - 3}\\
{ - 3}&{ - 4}&{ - 3}&{10}
\end{array}} \right].\]
The Coates graph of $\tilde{D}_{3,S}$ is shown in Figure 16.
\begin{figure}
\centering
\includegraphics[width=1.8in]{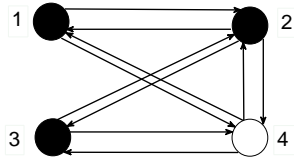}
\caption{The Coates graph of $\tilde{D}_{3,S}$}
\end{figure}
Let us take a look at the number of directed trees with reference node 1.
\[{\cal A}(G) = \left[ {\begin{array}{*{20}{c}}
0&1&0&1\\
1&0&1&1\\
0&1&0&1\\
1&1&1&0
\end{array}} \right],
{\cal L}(G) = \left[ {\begin{array}{*{20}{c}}
0&-1&2&-1\\
-1&3&-1&-1\\
0&-1&2&-1\\
-1&-1&-1&3
\end{array}} \right].\]
So the nonzero eigenvalues of ${\cal L}(G)$: ${\lambda _{_1}}=2$, ${\lambda _{_2}}= {\lambda _3}= 4$. Therefore, all directed trees with reference node 1 is
\[N = \frac{1}{4} \times 2 \times 4 \times 4 = 8.\]
Those directed trees with reference node 1 is shown in Figure 17-24.
\begin{figure}
\begin{minipage}[t]{0.5\linewidth}
\centering
\includegraphics[width=1.0in]{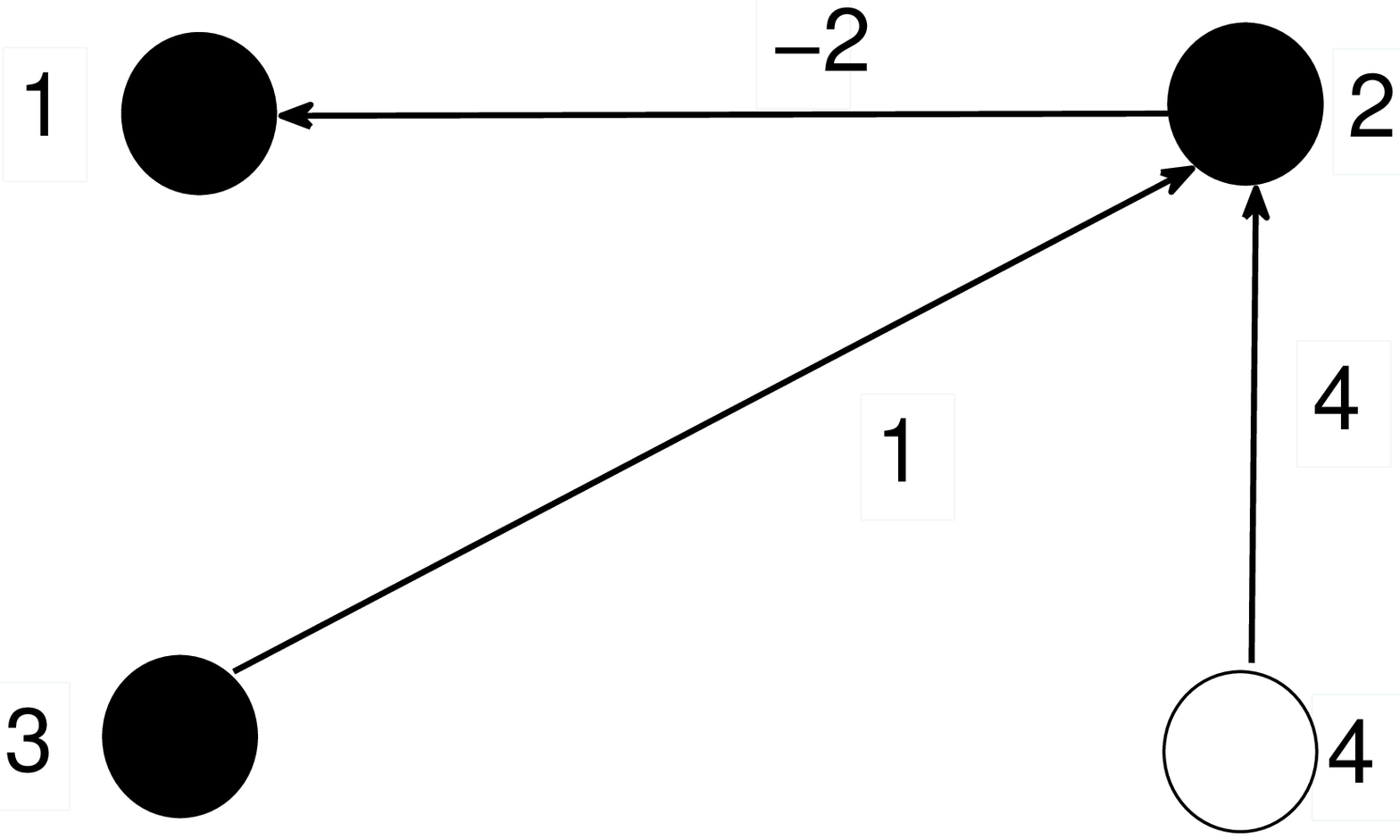}
\caption{$(1)\times(4)\times(-2)=-8 $}
\end{minipage}
\begin{minipage}[t]{0.5\linewidth}
\centering
\includegraphics[width=1.0in]{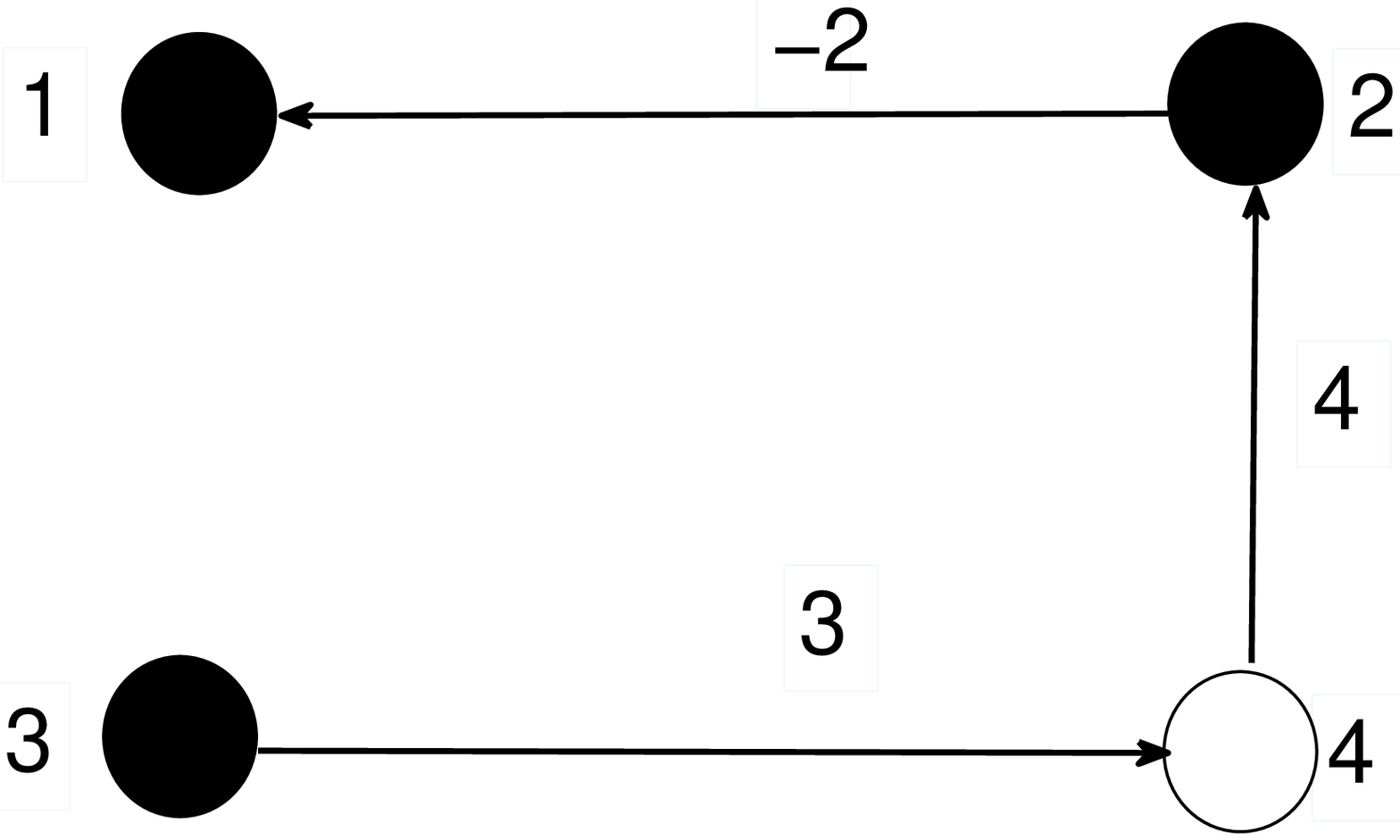}
\caption{$(3)\times(4)\times(-2)=-24$ }
\end{minipage}
\end{figure}

\begin{figure}
\begin{minipage}[t]{0.5\linewidth}
\centering
\includegraphics[width=1.0in]{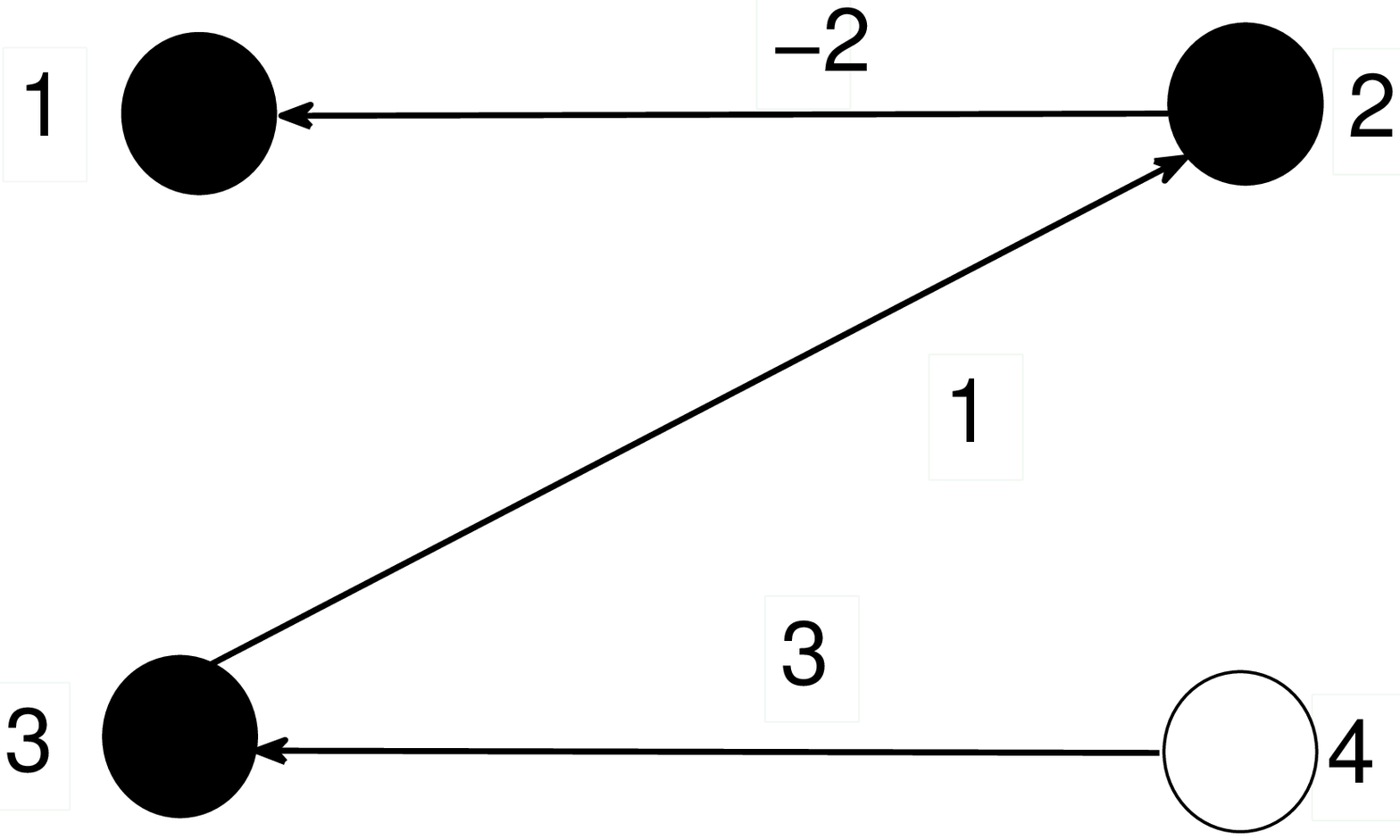}
\caption{$(3)\times(1)\times(-2)=-6$}
\end{minipage}
\begin{minipage}[t]{0.5\linewidth}
\centering
\includegraphics[width=1.0in]{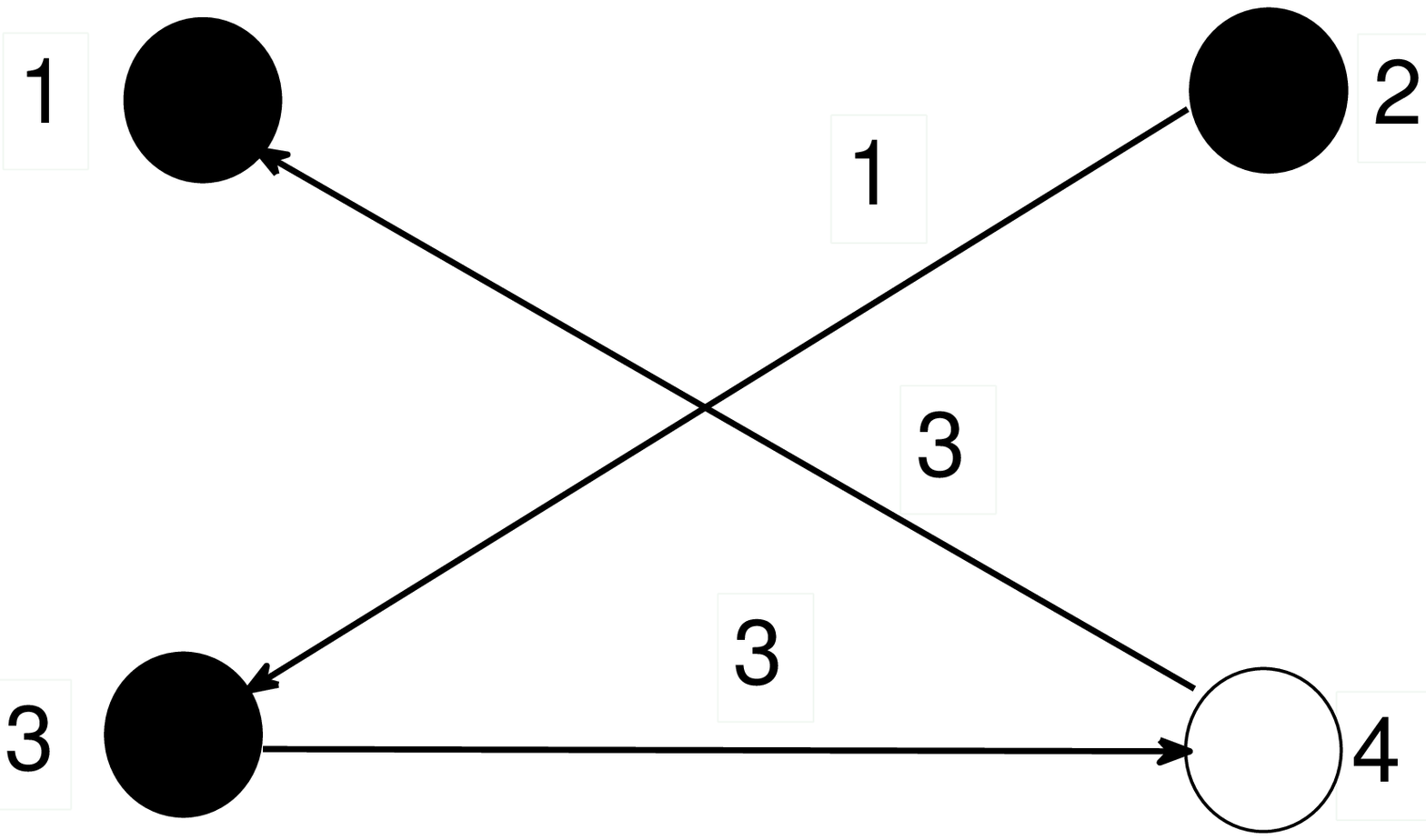}
\caption{$(1)\times(3)\times(3)=9$ }
\end{minipage}
\end{figure}

\begin{figure}
\begin{minipage}[t]{0.5\linewidth}
\centering
\includegraphics[width=1.0in]{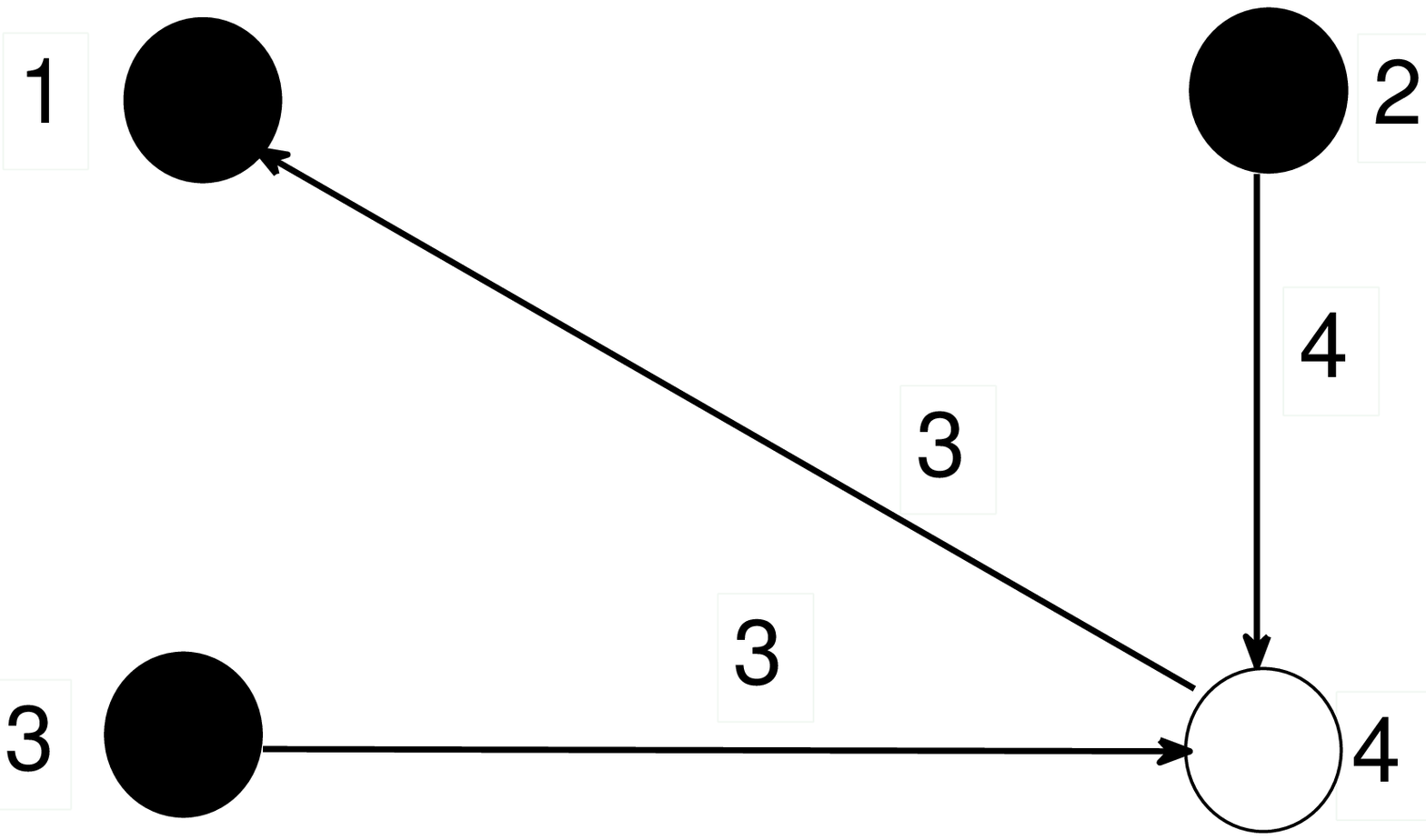}
\caption{$(3)\times(4)\times(3)=36$}
\end{minipage}
\begin{minipage}[t]{0.5\linewidth}
\centering
\includegraphics[width=1.0in]{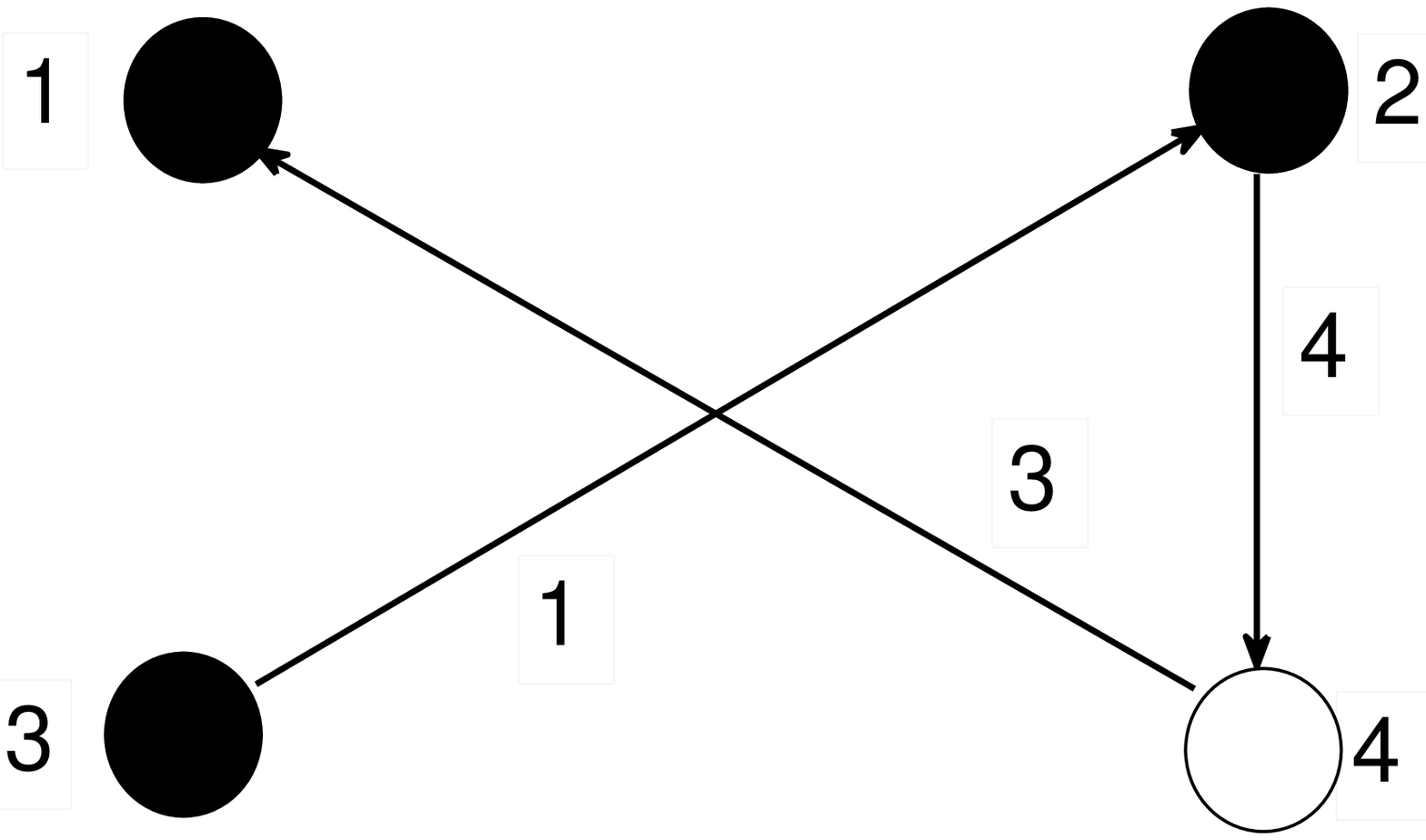}
\caption{$(1)\times(4)\times(3)=12$ }
\end{minipage}
\end{figure}

\begin{figure}
\begin{minipage}[t]{0.5\linewidth}
\centering
\includegraphics[width=1.0in]{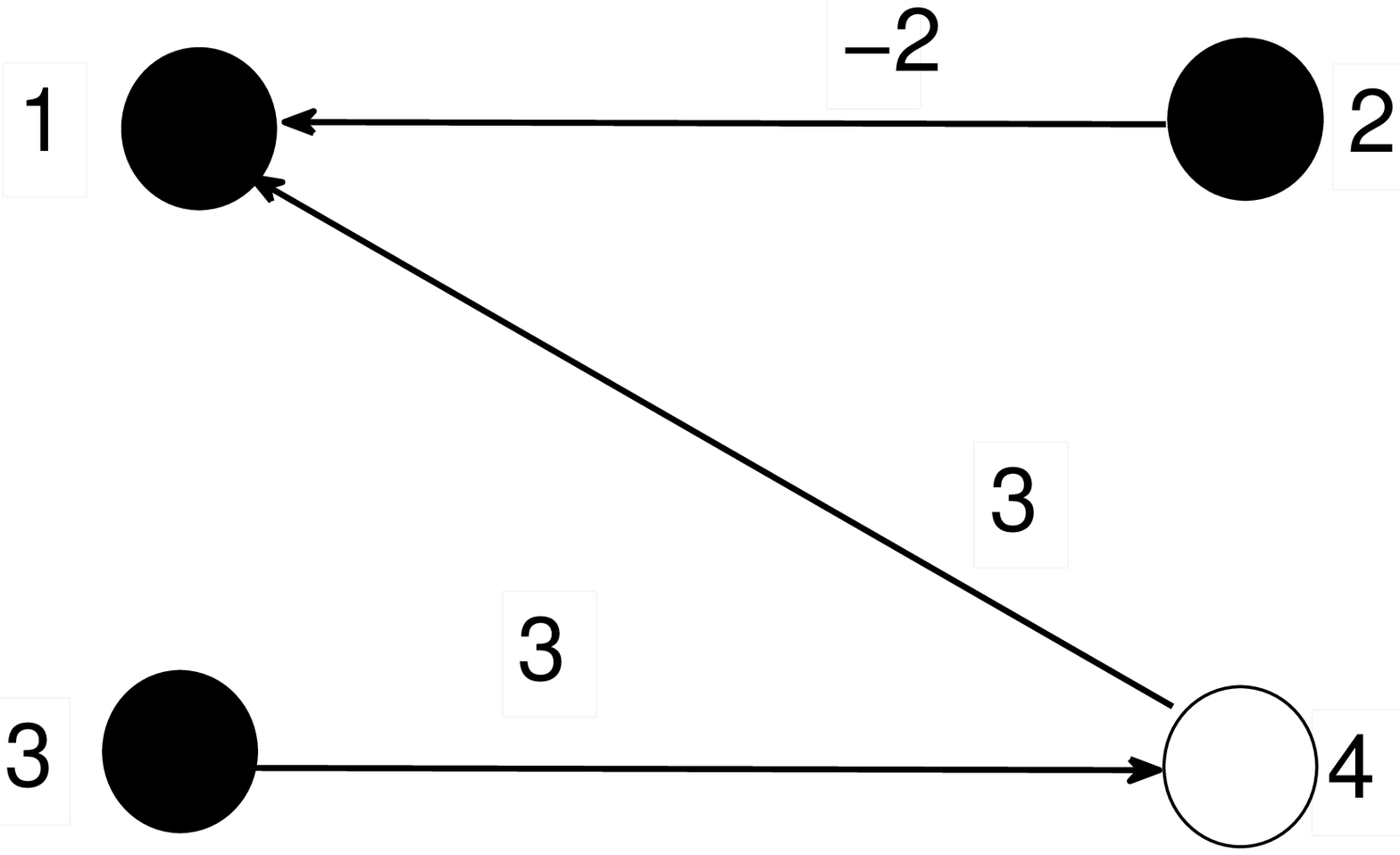}
\caption{$(3)\times(3)\times(-2)=-18$}
\end{minipage}
\begin{minipage}[t]{0.5\linewidth}
\centering
\includegraphics[width=1.0in]{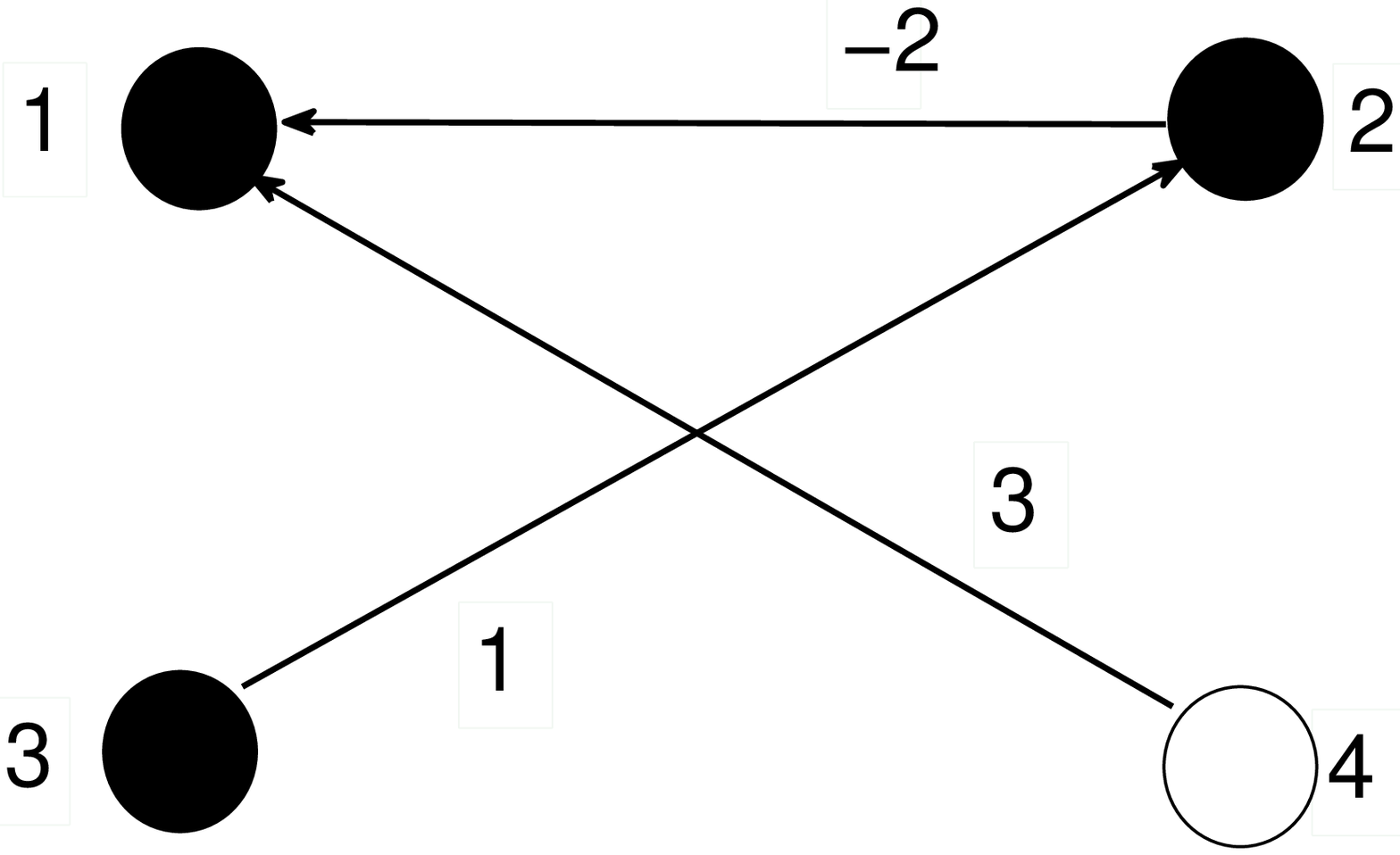}
\caption{$(1)\times(-2)\times(3)=-6$ }
\end{minipage}
\end{figure}
So
\[\det({D_{3,S}})=-6+(-8)+(-24)+(-6)+(-18)+(9)+(36)+(12)= - 5.\]

The $\Phi_S$ obtained by the PRL algorithm contains 13 cases, while the topology algorithm obtains the directed tree with reference node 1 have only 8 cases, which is much smaller than PRL algorithm. And even in some cases, there is only one case for some matrices. For example,
\[{\tilde{A}} = \left[ {\begin{array}{*{20}{c}}
1&2&0&{ - 3}\\
{ - 1}&0&1&0\\
0&0&2&{ - 2}\\
0&{ - 2}&{ - 3}&5
\end{array}} \right]\]
Matrix $\tilde{A}$ is an equicofactor matrix. The accompanying digraph of $A$ is shown in Figure 25. A directed tree with reference node 1 has only one, as shown in Figure 26.

\section{Conclusions}
In this paper, we clarify the relevant research results, obtain a more general representation, and provide a relevant algorithm framework. In terms of algorithms, the PRL algorithm is only applicable to symmetric matrices, however the topology algorithm can be used in any form of square matrix including asymmetric square matrices. Of course, our results can also be applied to reveal topological stability criteria for collective dynamics in complex networks, e.g., \cite{1,2,22}, which will continue to be studied.

\begin{figure}
\centering
\includegraphics[width=1.8in]{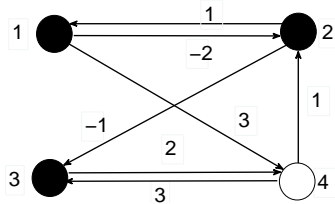}
\caption{The accompanying digraph of $\tilde{A}$}
\end{figure}

\begin{figure}[htp]
\centering
\includegraphics[width=1.8in]{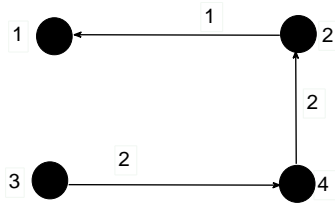}
\caption{The directed tree with reference node 1.}
\end{figure}


\end{document}